\documentclass[11pt, twoside]{article}
\usepackage{latexsym}
\usepackage{amsmath}
\usepackage{amssymb}
\usepackage[all]{xy}
\usepackage{amsfonts}
\usepackage{verbatim}
\usepackage{amsthm}
\usepackage{bm}
\usepackage{mathrsfs}
\usepackage{epsfig}
\usepackage{xy}
\usepackage{array}
\usepackage{stmaryrd}
\usepackage{graphicx,color}
\usepackage{xcolor}
\usepackage{tikz}
\usetikzlibrary{arrows,calc}
\usepackage{etex}
\usepackage{mathdots}
\usepackage{float}
\usepackage{graphics}
\usepackage{pdflscape}
\usepackage{CJK}
\usepackage{anysize,hyperref}
\input xypic
\xyoption{all}
\usepackage{extarrows}
\usepackage[perpage,symbol]{footmisc}
\usepackage{setspace}
\def\E{\mathbb{E}}
\def\B{\mathcal{B}}
\def\R{\mathcal{R}}
\def\Q{\mathcal{Q}}
\def\W{\mathcal{W}}
\def\V{\mathcal{V}}
\def\C{\mathcal{C}}
\def\s{\mathfrak{s}}

\def\op{^\mathrm{op}}
\def\Ab{\mathit{Ab}}
\def\del{\delta}
\def\dr{\ar@{->}[r]}

\def\Ext{\mbox{Ext}}
\def\Hom{\mbox{Hom}}

\begin{document}
\baselineskip=15pt
\title{\Large{\bf A Hovey triple arising from  two cotorsion pairs$^\bigstar$\footnotetext{\hspace{-1em}$^\bigstar$Panyue Zhou is supported by the National Natural Science Foundation of China (Grants No. 11901190 and 11671221) and the Hunan Provincial Natural Science Foundation of China (Grants No. 2018JJ3205).}}}
\medskip
\author{Panyue Zhou}

\date{}

\maketitle
\def\blue{\color{blue}}
\def\red{\color{red}}

\newtheorem{theorem}{Theorem}[section]
\newtheorem{lemma}[theorem]{Lemma}
\newtheorem{corollary}[theorem]{Corollary}
\newtheorem{proposition}[theorem]{Proposition}
\newtheorem{conjecture}{Conjecture}
\theoremstyle{definition}
\newtheorem{definition}[theorem]{Definition}
\newtheorem{question}[theorem]{Question}
\newtheorem{remark}[theorem]{Remark}
\newtheorem{remark*}[]{Remark}
\newtheorem{example}[theorem]{Example}
\newtheorem{example*}[]{Example}
\newtheorem{condition}[theorem]{Condition}
\newtheorem{condition*}[]{Condition}
\newtheorem{construction}[theorem]{Construction}
\newtheorem{construction*}[]{Construction}

\newtheorem{assumption}[theorem]{Assumption}
\newtheorem{assumption*}[]{Assumption}

\baselineskip=17pt
\parindent=0.5cm

\begin{abstract}
\baselineskip=16pt
Assume that $(\mathcal{C}, \mathbb{E}, \mathfrak{s})$ is an extriangulated category satisfying {\rm Condition (WIC)}.
Let $(\Q, \widetilde{\R})$ and $(\widetilde{\Q}, \R)$ be two hereditary cotorsion pairs satisfying conditions $\widetilde{\R} \subseteq \R$, $\widetilde{\Q}\subseteq \Q$ and $\widetilde{\Q}\cap \R = \Q \cap \widetilde{R}$. Then there exists a unique thick class $\W$ for which $(\Q,\W,\R)$ is a Hovey triple. As an application, this result generalizes the work by Gillespie in an exact case. Moreover, it highlights new phenomena when it applied to triangulated categories.\\[0.5cm]
\medskip
\end{abstract}

\pagestyle{myheadings}
\markboth{\rightline {\scriptsize   Panyue Zhou}}
         {\leftline{\scriptsize  A Hovey triple arising from  two cotorsion pairs}}

\section{Introduction}

Let $\C$ be an abelian category.
Hovey \cite{H1} showed that there exists a one-to-one correspondence between abelian model structures on $\C$ and triples $(\Q,\W,\R)$
having $\W$ is thick and  admitting two complete cotorsion pairs $(\Q,\W\cap\R)$ and $(\Q\cap\W,\R)$. In this case, the class $\Q$ is precisely the class of cofibrant objects, $\R$ is the class of fibrant objects, and $\W$ is the class of trivial objects in the abelian model structure on $\C$.
Given such a triple $(\Q,\W,\R)$, called a Hovey triple, denote the associated cotorsion pairs by $(\Q,\widetilde{\R})=(\Q,\W\cap\R)$ and $(\widetilde{\Q},\R)=(\Q\cap\W,\R)$. Then one has (1) $\widetilde{\R}\subseteq\R, \widetilde{Q}\subseteq\Q$ and (2) $\widetilde{Q}\cap\R=\Q\cap\widetilde{\R}$.  It holds in the
general setting of when $\C$ is a weakly idempotent complete exact category.
Gillespie \cite{G1} proved that Hovey's correspondence carries over to this
setting.

Gillespie \cite{G2} proved that there exists a converse when the cotorsion pairs above are hereditary. That is, if $(\Q,\widetilde{\R})$ and $(\widetilde{\Q},\R)$ are complete hereditary cotorsion pairs satisfying conditions (1) and (2) above, then there is a unique thick class $\W$ for which $(\Q,\W,\R)$ is a Hovey triple.

Extriangulated categories were recently introduced by Nakaoka and Palu \cite{NP} by
extracting those properties of $\Ext^1$ on exact categories (which is itself a generalisation of the concept
of a module category and an abelian category) and on triangulated categories
that seem relevant from the point of view of cotorsion pairs. In particular, exact
categories and triangulated categories are extriangulated categories. There are a lot of examples of extriangulated categories which are neither exact triangulated categories nor triangulated categories, see \cite{NP,ZZ}.
Hence, many results hold on exact categories and triangulated categories can be unified in the same framework.

Motivated by Hovey's correspondence, we introduce the notion of extriangulated model structures, which are a generalization of abelian model structures.

\begin{definition}
Let $(\C,\E,\s)$ be an extriangulated category, and $\Q,\W,\R$ three classes of objects in $\C$.
\begin{itemize}
\item[(1)] $\W$ is called \emph{thick} if it is closed under direct summands and satisfies
the $2$-out-of-$3$ property: for any $\E$-triangle $A\overset{}{\longrightarrow}B\overset{}{\longrightarrow}C\overset{}{\dashrightarrow}$
 in $\C$ with two terms in $\W$, the third term belongs to $\W$ as well.

\item[(2)] A triple $(\Q,\W,\R)$ is called \emph{extriangulated model structures} if $(\Q,\W\cap\R)$ and $(\Q\cap\W,\R)$ are two cotorsion pairs
and $\W$ is thick.  In this case, a triple $(\Q,\W,\R)$ is also called \emph{Hovey}.
\end{itemize}
\end{definition}
In this article we will generalise Gillespie's results \cite{G2} into the extriangulated case. Moreover,
our proof is not far from the usual exact case. We hope that our work will motivate and
further study on exangulated categories.

Our main result is the following.

\begin{theorem}\emph{(See Theorem \ref{main} for more details)}
Assume that $(\mathcal{C}, \mathbb{E}, \mathfrak{s})$ is an extriangulated category satisfying {\rm Condition (WIC)}.
Let $(\Q, \widetilde{\R})$ and $(\widetilde{\Q}, \R)$ be two hereditary cotorsion pairs satisfying conditions $\widetilde{\R} \subseteq \R$, $\widetilde{\Q}\subseteq \Q$ and $\widetilde{\Q}\cap \R = \Q \cap \widetilde{R}$. Then there exists a unique thick class $\W$ for which $(\Q,\W,\R)$ is a Hovey triple. Moreover, this thick class $\W$ can be described in the two following ways:
\begin{align*}
   \W  &= \{\, X \in \C \, | \,  \text{there exists an $\E$-triangle} \, X \to R \to Q\dashrightarrow \, \text{ with} \, R \in \widetilde{\R } \, , Q \in \widetilde{\Q} \,\} \\
           &= \{\, X \in \C \, | \,  \text{there exists an $\E$-triangle } \, R' \to Q' \to X \dashrightarrow\, \text{ with} \, R' \in \widetilde{\R} \, , Q' \in \widetilde{\Q} \,\}.
          \end{align*}
\end{theorem}

Our main result generalizes Gillespie's results on an exact category and
is new for a triangulated category case.

This article is organised as follows:
In Section 2, we review some basic concepts and results concerning extriangulated categories.
In Section 3, we show our main result.

\section{Preliminaries}
Let us briefly recall the definition and basic properties of extriangulated categories from \cite{NP}.

Let $\C$ be an additive category. Suppose that $\C$ is equipped with an additive bifunctor $$\E\colon\C\op\times\C\to\Ab,$$
where $\Ab$ is the category of abelian groups. For any pair of objects $A,C\in\C$, an element $\delta\in\E(C,A)$ is called an {\it $\E$-extension}. Thus formally, an $\E$-extension is a triplet $(A,\delta,C)$.
Let $(A,\del,C)$ be an $\E$-extension. Since $\E$ is a bifunctor, for any $a\in\C(A,A')$ and $c\in\C(C',C)$, we have $\E$-extensions
$$ \E(C,a)(\del)\in\E(C,A')\ \ \text{and}\ \ \ \E(c,A)(\del)\in\E(C',A). $$
We abbreviately denote them by $a_\ast\del$ and $c^\ast\del$.
For any $A,C\in\C$, the zero element $0\in\E(C,A)$ is called the \emph{spilt $\E$-extension}.

\begin{definition}{\cite[Definition 2.3]{NP}}
Let $(A,\del,C),(A',\del',C')$ be any pair of $\E$-extensions. A {\it morphism} $$(a,c)\colon(A,\del,C)\to(A',\del',C')$$ of $\E$-extensions is a pair of morphisms $a\in\C(A,A')$ and $c\in\C(C,C')$ in $\C$, satisfying the equality
$ a_\ast\del=c^\ast\del'. $
Simply we denote it as $(a,c)\colon\del\to\del'$.
\end{definition}

\begin{definition}{\cite[Definition 2.6]{NP}}
Let $\delta=(A,\delta,C),\delta^{\prime}=(A^{\prime},\delta^{\prime},C^{\prime})$ be any pair of $\mathbb{E}$-extensions. Let
\[ C\overset{\iota_C}{\longrightarrow}C\oplus C^{\prime}\overset{\iota_{C^{\prime}}}{\longleftarrow}C^{\prime} \]
and
\[ A\overset{p_A}{\longleftarrow}A\oplus A^{\prime}\overset{p_{A^{\prime}}}{\longrightarrow}A^{\prime} \]
be coproduct and product in $\B$, respectively. Since $\mathbb{E}$ is biadditive, we have a natural isomorphism
\[ \mathbb{E}(C\oplus C^{\prime},A\oplus A^{\prime})\cong \mathbb{E}(C,A)\oplus\mathbb{E}(C,A^{\prime})\oplus\mathbb{E}(C^{\prime},A)\oplus\mathbb{E}(C^{\prime},A^{\prime}). \]

Let $\delta\oplus\delta^{\prime}\in\mathbb{E}(C\oplus C^{\prime},A\oplus A^{\prime})$ be the element corresponding to $(\delta,0,0,\delta^{\prime})$ through the above isomorphism. This is the unique element which satisfies
\begin{eqnarray*}
\mathbb{E}(\iota_C,p_A)(\delta\oplus\delta^{\prime})=\delta&,&\mathbb{E}(\iota_C,p_{A^{\prime}})(\delta\oplus\delta^{\prime})=0,\\
\mathbb{E}(\iota_{C^{\prime}},p_A)(\delta\oplus\delta^{\prime})=0&,&\mathbb{E}(\iota_{C^{\prime}},p_{A^{\prime}})(\delta\oplus\delta^{\prime})=\delta^{\prime}.
\end{eqnarray*}
\end{definition}

\medskip

Let $A,C\in\C$ be any pair of objects. Sequences of morphisms in $\C$
$$\xymatrix@C=0.7cm{A\ar[r]^{x} & B \ar[r]^{y} & C}\ \ \text{and}\ \ \ \xymatrix@C=0.7cm{A\ar[r]^{x'} & B' \ar[r]^{y'} & C}$$
are said to be {\it equivalent} if there exists an isomorphism $b\in\C(B,B')$ which makes the following diagram commutative.
$$\xymatrix{
A \ar[r]^x \ar@{=}[d] & B\ar[r]^y \ar[d]_{\simeq}^{b} & C\ar@{=}[d]&\\
A\ar[r]^{x'} & B' \ar[r]^{y'} & C &}$$

We denote the equivalence class of $\xymatrix@C=0.7cm{A\ar[r]^{x} & B \ar[r]^{y} & C}$ by $[\xymatrix@C=0.7cm{A\ar[r]^{x} & B \ar[r]^{y} & C}]$.

For any $A,C\in\C$, we denote $0=[
A \xrightarrow{\binom{1}{0}} A\oplus C \xrightarrow{(0,1)} C ]$.

For any two equivalence classes, we denote
$$ [A\overset{x}{\longrightarrow}B\overset{y}{\longrightarrow}C]\oplus [A^{\prime}\overset{x^{\prime}}{\longrightarrow}B^{\prime}\overset{y^{\prime}}{\longrightarrow}C^{\prime}]=[A\oplus A^{\prime}\overset{x\oplus x^{\prime}}{\longrightarrow}B\oplus B^{\prime}\overset{y\oplus y^{\prime}}{\longrightarrow}C\oplus C^{\prime}]. $$

\begin{definition}{\cite[Definition 2.9]{NP}}
Let $\s$ be a correspondence which associates an equivalence class $\s(\del)=[\xymatrix@C=0.7cm{A\ar[r]^{x} & B \ar[r]^{y} & C}]$ to any $\E$-extension $\del\in\E(C,A)$. This $\s$ is called a {\it realization} of $\E$, if it satisfies the following condition:
\begin{itemize}
\item Let $\del\in\E(C,A)$ and $\del'\in\E(C',A')$ be any pair of $\E$-extensions, with $$\s(\del)=[\xymatrix@C=0.7cm{A\ar[r]^{x} & B \ar[r]^{y} & C}],\ \ \ \s(\del')=[\xymatrix@C=0.7cm{A'\ar[r]^{x'} & B'\ar[r]^{y'} & C'}].$$
Then, for any morphism $(a,c)\colon\del\to\del'$, there exists $b\in\C(B,B')$ which makes the following diagram commutative.
$$\xymatrix{
A \ar[r]^x \ar[d]^a & B\ar[r]^y \ar[d]^{b} & C\ar[d]^c&\\
A'\ar[r]^{x'} & B' \ar[r]^{y'} & C' &}$$

\end{itemize}
In the above situation, we say that the triplet $(a,b,c)$ realizes $(a,c)$.
\end{definition}

\begin{definition}{\cite[Definition 2.10]{NP}}
A realization $\s$ of $\E$ is called \emph{additive} if it satisfies the following conditions.
\begin{itemize}
\item[(1)] For any $A,C\in\C$, the split $\E$-extension $0\in\E(C,A)$ satisfies $\s(0)=0$.
\item[(2)] For any pair of $\E$-extensions $\delta\in\E(C,A)$ and $\delta'\in\E(C',A')$,
$$\s(\delta\oplus\delta')=\s(\delta)\oplus\s(\delta')$$
holds.
\end{itemize}
\end{definition}

\begin{definition}{\cite[Definition 2.12]{NP}}
A triplet $(\C,\E,\s)$ is called an \emph{externally triangulated category} (or \emph{extriangulated category} for short) if it satisfies the following conditions:
\begin{itemize}
\item[{\rm (ET1)}] $\E\colon\C\op\times\C\to\Ab$ is a biadditive functor.
\item[{\rm (ET2)}] $\s$ is an additive realization of $\E$.
\item[{\rm (ET3)}] Let $\del\in\E(C,A)$ and $\del'\in\E(C',A')$ be any pair of $\E$-extensions, realized as
$$ \s(\del)=[\xymatrix@C=0.7cm{A\ar[r]^{x} & B \ar[r]^{y} & C}],\ \ \s(\del')=[\xymatrix@C=0.7cm{A'\ar[r]^{x'} & B' \ar[r]^{y'} & C'}]. $$
For any commutative square
$$\xymatrix{
A \ar[r]^x \ar[d]^a & B\ar[r]^y \ar[d]^{b} & C&\\
A'\ar[r]^{x'} & B' \ar[r]^{y'} & C' &}$$
in $\C$, there exists a morphism $(a,c)\colon\del\to\del'$ satisfying $cy=y'b$.
\item[{\rm (ET3)$\op$}] Dual of (ET3).
\item[{\rm (ET4)}] Let $(A,\del,D)$ and $(B,\del',F)$ be $\E$-extensions realized by
$$\xymatrix@C=0.7cm{A\ar[r]^{f} & B \ar[r]^{f'} & D}\ \ \text{and}\ \ \ \xymatrix@C=0.7cm{B\ar[r]^{g} & C \ar[r]^{g'} & F}$$
respectively. Then there exist an object $E\in\C$, a commutative diagram
$$\xymatrix{A\ar[r]^{f}\ar@{=}[d]&B\ar[r]^{f'}\ar[d]^{g}&D\ar[d]^{d}\\
A\ar[r]^{h}&C\ar[d]^{g'}\ar[r]^{h'}&E\ar[d]^{e}\\
&F\ar@{=}[r]&F}$$
in $\C$, and an $\E$-extension $\del^{''}\in\E(E,A)$ realized by $\xymatrix@C=0.7cm{A\ar[r]^{h} & C \ar[r]^{h'} & E},$ which satisfy the following compatibilities.
\begin{itemize}
\item[{\rm (i)}] $\xymatrix@C=0.7cm{D\ar[r]^{d} & E \ar[r]^{e} & F}$  realizes $f'_{\ast}\del'$,
\item[{\rm (ii)}] $d^\ast\del''=\del$,

\item[{\rm (iii)}] $f_{\ast}\del''=e^{\ast}\del'$.
\end{itemize}

\item[{\rm (ET4)$\op$}]  Dual of (ET4).
\end{itemize}
\end{definition}

We will use the following terminology.
\begin{definition}{\cite{NP}}
Let $(\C,\E,\s)$ be an extriangulated category.
\begin{itemize}
\item[(1)] A sequence $A\xrightarrow{~x~}B\xrightarrow{~y~}C$ is called a {\it conflation} if it realizes some $\E$-extension $\del\in\E(C,A)$.
    In this case, $x$ is called an {\it inflation} and $y$ is called a {\it deflation}.

\item[(2)] If a conflation  $A\xrightarrow{~x~}B\xrightarrow{~y~}C$ realizes $\delta\in\mathbb{E}(C,A)$, we call the pair $( A\xrightarrow{~x~}B\xrightarrow{~y~}C,\delta)$ an {\it $\E$-triangle}, and write it in the following way.
$$A\overset{x}{\longrightarrow}B\overset{y}{\longrightarrow}C\overset{\delta}{\dashrightarrow}$$
We usually do not write this $``\delta"$ if it is not used in the argument.

\item[(3)] Let $A\overset{x}{\longrightarrow}B\overset{y}{\longrightarrow}C\overset{\delta}{\dashrightarrow}$ and $A^{\prime}\overset{x^{\prime}}{\longrightarrow}B^{\prime}\overset{y^{\prime}}{\longrightarrow}C^{\prime}\overset{\delta^{\prime}}{\dashrightarrow}$ be any pair of $\E$-triangles. If a triplet $(a,b,c)$ realizes $(a,c)\colon\delta\to\delta^{\prime}$, then we write it as
$$\xymatrix{
A \ar[r]^x \ar[d]^a & B\ar[r]^y \ar[d]^{b} & C\ar@{-->}[r]^{\del}\ar[d]^c&\\
A'\ar[r]^{x'} & B' \ar[r]^{y'} & C'\ar@{-->}[r]^{\del'} &}$$
and call $(a,b,c)$ a {\it morphism of $\E$-triangles}.
\end{itemize}
\end{definition}

Assume that $(\C, \E, \s)$ is an extriangulated category. By Yoneda's lemma, any $\E$-extension
$\del\in\E(C,A)$ induces natural transformations
$$\del_{\sharp}\colon\C(-,C)\Rightarrow\E(-,A)\ \ \textrm{and} \ \ \del^{\sharp}\colon \C(A,-)\Rightarrow\E(C,-).$$
For any $X\in\C$, these $(\del_{\sharp})_X $ and $\del_{X}^{\sharp}$ are given as follows:

(1) $(\del_{\sharp})_X\colon\C(X,C)\to \E(X,A); f\mapsto f^{\ast}\del$.

(2) $\del_{X}^{\sharp}\colon \C(A,X)\to \E(C,X); g\mapsto g_{\ast}\del$.

When there is no danger of confusion, we will sometimes instead of
$$\Hom_\C(X,A)\xrightarrow{~\Hom_\C(X,\ f)~}\Hom_\C(X,B)$$
 write the simplified form:
$$\C(X,A)\xrightarrow{~\C(X,\ f)~}\C(X,B).$$

\begin{lemma}\label{lem}
Let $(\C, \E, \s)$ be an extriangulated category, and $$\xymatrix{A\ar[r]^{x}&B\ar[r]^{y}&C\ar@{-->}[r]^{\delta}&}$$
an $\E$-triangle. Then we have the following long exact sequence:
$$\C(-, A)\xrightarrow{\C(-,x)}\C(-, B)\xrightarrow{\C(-,y)}\C(-, C)\xrightarrow{\delta^{\sharp}_-}
\E(-, A)\xrightarrow{\E(-,x)}\E(-, B)\xrightarrow{\E(-,y)}\E(-, C);
$$
$$\C(C,-)\xrightarrow{\C(y,-)}\C(B,-)\xrightarrow{\C(x,-)}\C(A,-)\xrightarrow{\delta_{\sharp}^-}
\E(C,-)\xrightarrow{\E(y,-)}\E(B,-)\xrightarrow{\E(x,-)}\E(A,-)
.$$
\end{lemma}

\proof This follows from Proposition 3.3 and Proposition 3.11 in \cite{NP}. \qed
\medskip

We recall the following result, which (also the dual of it) will be used later.

\begin{lemma}\emph{\cite[Proposition 1.20]{LN}}
\label{y2} Assume that $(\C, \E, \s)$ is an extriangulated category.
Let $A\overset{x}{\longrightarrow}B\overset{y}{\longrightarrow}C\overset{\delta}{\dashrightarrow}$ be any $\E$-triangle, and $f\colon A\rightarrow D$ be any morphism in $\C$. Then there exists a morphism $g$ which gives a morphism of $\E$-triangles
$$\xymatrix{
A \ar[r]^{x} \ar[d]_f &B \ar[r]^{y} \ar[d]^g &C \ar@{=}[d]\ar@{-->}[r]^{\delta}&\\
D \ar[r]^{d} &F \ar[r]^{e} &C\ar@{-->}[r]^{f_{\ast}\delta}&
}
$$
and moreover, the sequence $A\xrightarrow{\binom{-x}{f}}B\oplus D\xrightarrow{(g,\ d)}F\overset{f_{\ast}\delta}{\dashrightarrow}$ becomes an $\E$-triangle.
\end{lemma}

\begin{lemma}\label{y3}
Assume that $(\C, \E, \s)$ is an extriangulated category. Let
$$\xymatrix{A\ar[r]^f\ar[d]^{a}&B\ar[r]^g\ar[d]^{b}&C\ar@{-->}[r]^{\delta}\ar@{=}[d]&\\
D\ar[r]^{u}&F\ar[r]^{v}&C\ar@{-->}[r]^{\eta}&}$$
be a morphism of $\E$-triangles.
Then there exists an isomorphism $w\colon F\to F$ which makes
$$A\xrightarrow{\binom{-f}{a}}B\oplus D\xrightarrow{(b,\ u)}F\overset{\theta}{\dashrightarrow}$$
an $\E$-triangle, where $\theta=(w^{-1})^\ast v^{\ast}\del$.
\end{lemma}

\proof By Lemma \ref{y2}, we have the following commutative diagram of $\E$-triangles
$$\xymatrix{A\ar[r]^f\ar[d]^{a}&B\ar[r]^g\ar[d]^{b'}&C\ar@{-->}[r]^{\delta}\ar@{=}[d]&\\
D\ar[r]^{x}&F'\ar[r]^{y}&C\ar@{-->}[r]^{a_{\ast}\del}&}$$
and moreover,
$$A\xrightarrow{\binom{-f}{a}}B\oplus D\xrightarrow{(b',\ x)}F'\overset{y^{\ast}\del}{\dashrightarrow}$$
is an $\E$-triangle.
By the definition of realization and \cite[Corollary 3.6]{NP}, there exists an isomorphism
$c\colon F'\to F$ which gives a morphism of $\E$-triangles.
$$\xymatrix{
D\ar[r]^{x}\ar@{=}[d]&F'\ar[r]^{y}\ar@{-->}[d]^{c}&C\ar@{=}[d]\ar@{-->}[r]^{a_{\ast}\del}&\\
D\ar[r]^{u}&F\ar[r]^{v}&C\ar@{-->}[r]^{\eta}&.}$$
Consider the isomorphism $c^{-1}\colon F\to F'$, by \cite[Proposition 3.7]{NP}, we obtain that
$$A\xrightarrow{\binom{-f}{a}}B\oplus D\xrightarrow{(cb',\ u)}F\overset{v^{\ast}\del}{\dashrightarrow}$$
is an $\E$-triangle.
Note that $(b,u)\binom{-f}{a}=0$, there exists a morphism $w\colon F\to F$ such that
$(b,u)=w(cb',u)$ and then $b=wcb'$ and $u=wu$. Since $(1-w)u=0$, there exists  a morphism $d\colon C\to F$ such that
$1-w=dv$. It follows that $$(1+dv)(1-dv)u=(1+dv)wu=(1+dv)u=u.$$

Observe that $v(b-cb')=vb-vcb'=g-yb'=0$, there exists a morphism $m\colon B\to D$ such that $b-cb'=um$ and then
$dvb=dv(um+cb')=dvcb'$. Thus we have $(1+dv)b=b+dvb=wcb'+dvcb'=(w+dv)cb'=cb'$ implies
$$(1+dv)(1-dv)cb'=(1+dv)wcb'=(1+dv)b=cb'.$$
Hence we have the following commutative diagram
$$\xymatrix@C=1.2cm{
A\ar[r]^{\binom{-f}{a}\quad}\ar@{=}[d]&B\oplus D\ar[r]^{\quad (cb',\ u)}\ar@{=}[d]&F\ar[d]^{(1+dv)(1-dv)}\ar@{-->}[r]^{v^{\ast}\del}&\\
A\ar[r]^{\binom{-f}{a}\quad}&B\oplus D\ar[r]^{\quad (cb',\ u)}&F\ar@{-->}[r]^{v^{\ast}\del}&.}$$
of $\E$-triangles. By \cite[Corollary 3.6]{NP}, we know that $(1+dv)(1-dv)$ is an isomorphism.  Note that $(1+dv)(1-dv)=(1-dv)(1+dv)$, we obtain that
$1-dv=w$ is an isomorphism. Consider the isomorphism $w^{-1}\colon F\to F$, by \cite[Proposition 3.7]{NP}, we get that
$$A\xrightarrow{\binom{-f}{a}}B\oplus D\xrightarrow{(b,\ u)}F\overset{\theta}{\dashrightarrow}$$
is an $\E$-triangle, where $\theta=(w^{-1})^\ast v^{\ast}\del$. \qed

\section{A Hovey triple arising from two cotorsion pairs}
In this section, we will prove our main result, before we need some preparations as follows.

The following conditions are from \cite[Condition 5.8]{NP}, analogous to the weak idempotent completeness.

\begin{condition}({\rm \textbf{Condition (WIC)}}) Let  $(\C, \E, \s)$ be an extriangulated category. Consider the following conditions.
\begin{itemize}
\item[(1)] Let $f\in\mathcal{C}(A, B),~ g\in\mathcal{C}(B, C)$ be any composable pair of morphisms. If $gf$ is an inflation, then so is $f$.

\item[(2)] Let $f\in\mathcal{C}(A, B), ~g\in\mathcal{C}(B, C)$ be any composable pair of morphisms. If $gf$ is a deflation, then so is $g$.
\end{itemize}
\end{condition}

\begin{remark}
(1) If $\C$ is an exact category, then Condition (WIC) is equivalent to $\C$ is
weakly idempotent complete, see \cite[Proposition 7.6]{B}.

(2) If $\C$ is a triangulated category, then Condition (WIC) is automatically satisfied.
\end{remark}

\begin{lemma}\label{y4}
Assume that $(\C, \E, \s)$ is an extriangulated category satisfying {\rm Condition (WIC)}. Let
$$\xymatrix{
A \ar[r]^x \ar@{=}[d] & B\ar[r]^y \ar[d]^{b} & C\ar@{-->}[r]^{\del}\ar[d]^c&\\
A\ar[r]^{x'} & B' \ar[r]^{y'} & C'\ar@{-->}[r]^{\del'} &}$$
be a morphism of $\E$-triangles. If  $c$ is a deflation, then $b$ is also a deflation.
\end{lemma}

\proof Since $c$ is a deflation, there exists an $\E$-triangle
$D\overset{d}{\longrightarrow}C\overset{c}{\longrightarrow}C'\overset{\eta}{\dashrightarrow}$
in $\C$. By \cite[Proposition 3.15]{NP}, we obtain the following commutative diagram made of $\E$-triangles.
$$\xymatrix{&D\ar@{=}\ar[d]^{e}\ar@{=}[r]&D\ar[d]^{d}&\\
A\ar@{=}[d]\ar[r]^{g} & M\ar[r]^h\ar[d]^{f} &
C\ar[d]^c\ar@{-->}[r]&\\
 A'\ar[r]^{x'} & B'\ar@{-->}[d]\ar[r]^{y'} & C'\ar@{-->}[d]^{\eta}\ar@{-->}[r]^{\delta'}&\\
 &&&}$$
It follows that $ch=y'f$ and then $(c,y')\binom{-h}{f}=0$.
By the dual of Lemma \ref{y3}, we have that
$$B\xrightarrow{\binom{-y}{b}}C\oplus B'\xrightarrow{(c,\ y')}C'\overset{\theta}{\dashrightarrow}$$
is an $\E$-triangle. So there exists a morphism $k\colon M\to B$ such that
$\binom{-y}{b}k=\binom{-h}{f}.$
In particular, we have $f=bk$.  By Condition (WIC) and $f$ is a deflation, this
$b$ becomes a deflation.  \qed

\begin{lemma}\label{y5}
Assume that $(\C, \E, \s)$ is an extriangulated category satisfying {\rm Condition (WIC)}. Let
$$\xymatrix{
A \ar[r]^x \ar[d]^a & B\ar[r]^y \ar[d]^{b} & C\ar@{-->}[r]^{\del}\ar[d]^c&\\
A'\ar[r]^{x'} & B' \ar[r]^{y'} & C'\ar@{-->}[r]^{\del'} &}$$
be a morphism of $\E$-triangles.
If $a$ and $c$ are deflations, then there exists a deflation $b'\colon B\to B'$ which gives the following
morphism of $\E$-triangles.
$$\xymatrix{
A \ar[r]^x \ar[d]^a & B\ar[r]^y \ar@{-->}[d]^{b'} & C\ar@{-->}[r]^{\del}\ar[d]^c&\\
A'\ar[r]^{x'} & B' \ar[r]^{y'} & C'\ar@{-->}[r]^{\del'} &}$$
\end{lemma}

\proof Since $a$ is a deflation, there exists an $\E$-triangle
$D\overset{d}{\longrightarrow}A\overset{a}{\longrightarrow}A'\overset{\eta}{\dashrightarrow}$
in $\C$. By (ET4),  we have the following commutative diagram
$$\xymatrix{D\ar[r]^d\ar@{=}[d] & A\ar[r]^a\ar[d]^x &A'\ar@{-->}[r]^{\eta}\ar[d]^u&\\
D\ar[r]^s & B\ar[r]^t\ar[d]^y & M\ar@{-->}[r]\ar[d]^v &\\
& C\ar@{-->}[d]^{\delta}\ar@{=}[r] & C\ar@{-->}[d]^{a_{\ast}\delta}\\&&}$$
of $\E$-triangles. By the dual of Lemma \ref{y2}, we have the following commutative diagram
$$\xymatrix{
A' \ar[r]^p \ar@{=}[d] & L\ar[r]^q \ar[d]^{g} & C\ar@{-->}[r]^{c^{\ast}\del'}\ar[d]^c&\\
A'\ar[r]^{x'} & B' \ar[r]^{y'} & C'\ar@{-->}[r]^{\del'} &}$$
of $\E$-triangles. By Lemma \ref{y4}, we have that $g$ is a deflation since $c$ is a deflation.
By the definition of realization and \cite[Corollary 3.6]{NP}, there exists an
isomorphism $h\colon M\to L$ which gives a morphism of $\E$-triangles.
$$\xymatrix{
A' \ar[r]^u \ar@{=}[d] & M\ar[r]^v \ar@{-->}[d]^{h} & C\ar@{-->}[r]^{a_{\ast}\del}\ar@{=}[d]&\\
A'\ar[r]^{p} & L \ar[r]^{q} & C\ar@{-->}[r]^{c^{\ast}\del'} &}$$
Hence we have the following commutative diagram
$$\xymatrix{
A \ar[r]^x \ar[d]^a & B\ar[r]^y \ar[d]^{ght} & C\ar@{-->}[r]^{\del}\ar[d]^c&\\
A'\ar[r]^{x'} & B' \ar[r]^{y'} & C'\ar@{-->}[r]^{\del'} &}$$
of $\E$-triangles. Since $t,h$ and $g$ are deflations,  we obtain that $ght$ is a deflation.   Put $b':=ght$, it is what we want.  \qed
\medskip

The following result is similar to the nine lemma.

\begin{lemma}\emph{\cite[Lemma 5.9]{NP}}
Assume that $(\C,\E,\s)$ is an extriangulated category satisfying {\rm Condition (WIC)}. Let
$$\xymatrix{
  K \ar[d]_{k}  & K' \ar[d]_{k'}  &  & \\
  A\ar@{}[dr]|{\circlearrowleft} \ar[d]_{a} \ar[r]^{x} & B \ar@{}[dr]|{\circlearrowleft}\ar[d]_{b} \ar[r]^{y} & C  \ar@{-->}[r]^{\delta} &  \\
  A' \ar@{-->}[d]_{\kappa} \ar[r]_{x'} & B' \ar@{-->}[d]_{\kappa'} \ar[r]^{y'} & C' \ar@{-->}[r]^{\delta'} &  \\
  &  & &  }
  $$
 be a diagram made of $\E$-triangles.
 Then for some $X\in\C$, we obtain  $\mathbb{E}$-triangles
 $$\xymatrix{K\ar[r]^m&K'\ar[r]^{n}&X\ar@{-->}[r]^{\nu}&}~\textrm{and}~\xymatrix{X\ar[r]^i&C\ar[r]^{c}&C'\ar@{-->}[r]^{\tau}&}$$ which make the following diagram commutative,
 $$ \xymatrix{
  K \ar@{}[dr]|{\circlearrowleft}\ar[d]_{k}\ar[r]^m  & K'\ar@{}[dr]|{\circlearrowleft} \ar[d]_{k'}\ar[r]^{n}  & X \ar[d]_{i}\ar@{-->}[r]^{\nu}  & \\
  A \ar@{}[dr]|{\circlearrowleft}\ar[d]_{a} \ar[r]^{x} & B\ar@{}[dr]|{\circlearrowleft} \ar[d]_{b} \ar[r]^{y} & C \ar[d]_{c} \ar@{-->}[r]^{\delta} &  \\
  A' \ar@{-->}[d]_{\kappa} \ar[r]^{x'} & B' \ar@{-->}[d]_{\kappa'} \ar[r]^{y'} & C' \ar@{-->}[d]_{\tau} \ar@{-->}[r]^{\delta'} &  \\
  &  & &  }
  $$
in which, those $(k, k', i), (a,b,c),(m,x,x')$ and $(n, y, y')$ are morphisms of $\mathbb{E}$-triangles.
\end{lemma}

\begin{lemma}\label{y8}
 Assume that $(\mathcal{C}, \mathbb{E}, \mathfrak{s})$ is an extriangulated category satisfying {\rm Condition (WIC)}. Let
$$
\xymatrix{
  K \ar[d]_{k}  & K' \ar[d]_{k'}  & K'' \ar[d]_{k''}  & \\
  A\ar@{}[dr]|{\circlearrowleft} \ar[d]_{a} \ar[r]^{x} & B \ar@{}[dr]|{\circlearrowleft}\ar[d]_{b} \ar[r]^{y} & C \ar[d]_{c'} \ar@{-->}[r]^{\delta} &  \\
  A' \ar@{-->}[d]_{\kappa} \ar[r]_{x'} & B' \ar@{-->}[d]_{\kappa'} \ar[r]^{y'} & C' \ar@{-->}[d]_{\kappa''} \ar@{-->}[r]^{\delta'} &  \\
  &  & &  }
  $$
 be a diagram made of $\mathbb{E}$-triangles.  If $y$ is a retraction, then there exists an $\mathbb{E}$-triangle $$\xymatrix{K\ar[r]^m&K'\ar[r]^{n'}&K''\ar@{-->}[r]^{\delta''}&}$$
 which make the following diagram commutative
 $$
\xymatrix{
  K\ar@{}[dr]|{\circlearrowleft} \ar[d]_{k}\ar[r]^m  & K' \ar@{}[dr]|{\circlearrowleft}\ar[d]_{k'}\ar[r]^{n'}  & K'' \ar[d]_{k''}\ar@{-->}[r]^{\delta''}  & \\
  A\ar@{}[dr]|{\circlearrowleft} \ar[d]_{a} \ar[r]^{x} & B \ar@{}[dr]|{\circlearrowleft}\ar[d]_{b} \ar[r]^{y} & C \ar[d]_{c'} \ar@{-->}[r]^{\delta} &  \\
  A' \ar@{-->}[d]_{\kappa} \ar[r]^{x'} & B' \ar@{-->}[d]_{\kappa'} \ar[r]^{y'} & C' \ar@{-->}[d]_{\kappa''} \ar@{-->}[r]^{\delta'} &  \\
  &  & &  }
  $$
in which, those $(k, k', k''), (a,b,c), (m, x, x')$ and $(n', y, y')$ are morphisms of $\mathbb{E}$-triangles.
\end{lemma}

\proof By Lemma \ref{y8}, there are $\mathbb{E}$-triangles
$$\xymatrix{K\ar[r]^m&K'\ar[r]^{n}&X\ar@{-->}[r]^{\nu}&}~\textrm{and}~
\xymatrix{X\ar[r]^i&C\ar[r]^{c}&C'\ar@{-->}[r]^{\tau}&}$$ which make the following diagram commutative
 $$ \xymatrix{
  K \ar@{}[dr]|{\circlearrowleft}\ar[d]_{k}\ar[r]^m  & K'\ar@{}[dr]|{\circlearrowleft} \ar[d]_{k'}\ar[r]^{n}  & X \ar[d]_{i}\ar@{-->}[r]^{\nu}  & \\
  A \ar@{}[dr]|{\circlearrowleft}\ar[d]_{a} \ar[r]^{x} & B\ar@{}[dr]|{\circlearrowleft} \ar[d]_{b} \ar[r]^{y} & C \ar[d]_{c} \ar@{-->}[r]^{\delta} &  \\
  A' \ar@{-->}[d]_{\kappa} \ar[r]^{x'} & B' \ar@{-->}[d]_{\kappa'} \ar[r]^{y'} & C' \ar@{-->}[d]_{\tau} \ar@{-->}[r]^{\delta'} &  \\
  &  & &  }
  $$
in which, those $(k, k', i), (a, b, c), (m, x, x')$ and $(n, y, y')$ are morphisms of $\mathbb{E}$-triangles.
Since $c'y=y'b=cy$ and $y$ is a retraction, we have $c'=c$. By (ET3)$^{\rm op}$ and \cite[Corollary 3.6]{NP}, we obtain a morphism of $\mathbb{E}$-triangles
$$\xymatrix{K''\ar[r]^{k''}\ar[d]^u&C\ar@{=}[d]\ar[r]^{c'}&C'\ar@{=}[d]\ar@{-->}[r]^{\kappa''}&\\
X\ar[r]^i&C\ar[r]^{c}&C'\ar@{-->}[r]^{\tau}&\\
}$$
where $u$ is an isomorphism. It is clear that $u_*\kappa''=\tau$. Put $n'=u^{-1}n$, by \cite[Proposition 3.7]{NP}, there exists an $\E$-triangle $$\xymatrix{K\ar[r]^{m\;}&K'\ar[r]^{n'}&K''\ar@{-->}[r]^{\delta''}&}$$
where $\delta''=u^*\nu$.
It is straightforward to verify that $k''n'=yk'$. Moreover, one can show that $$k_*\delta''=k_*u^*\nu=u^*k_*\nu=u^*i^*\delta=(iu)^*\delta=k^{''*}\delta,$$ $$n'_*\kappa'=u^{-1}_*n_*\kappa'=u^{-1}_*{y'}^{*}\tau={y'}^{*}u^{-1}\tau={y'}^{*}\kappa'',$$ which implies that $(k, k', k'')$ and $(n', y, y')$ are morphisms of $\mathbb{E}$-triangles.  \qed

\begin{definition}
Let $(\C,\E,\s)$ be an extriangulated category and $\W$ a class of objects in $\C$.
\begin{itemize}
\item[(1)] $\W$ is closed under cocones of deflations if for any $\E$-triangle $A\overset{x}{\longrightarrow}B\overset{y}{\longrightarrow}C\overset{\delta}{\dashrightarrow}$
which satisfies $B,C\in\W$, then $A\in\W$.
\item[(2)] $\W$ is closed under cones of inflations if for any $\E$-triangle $A\overset{x}{\longrightarrow}B\overset{y}{\longrightarrow}C\overset{\delta}{\dashrightarrow}$
which satisfies $A,B\in\W$, then $C\in\W$.
\item[(3)] $\W$ is closed under extensions if for any $\E$-triangle $A\overset{x}{\longrightarrow}B\overset{y}{\longrightarrow}C\overset{\delta}{\dashrightarrow}$
which satisfies $A,C\in\W$, then $B\in\W$.
\end{itemize}
\end{definition}

\begin{definition}\cite[Definition 4.1]{NP}
Assume that $(\C,\E,\s)$ is an extriangulated category.
Let $\Q$ and $\R$ be two classes of objects in $\C$. We call $(\Q,\R)$ a \emph{cotorsion pair} if it satisfies the following conditions:
\begin{itemize}
\setlength{\itemsep}{0.5pt}
\item[{\rm (a)}] $\E(\Q,\R)=0$.

\item[{\rm (b)}] For any object $C\in \C$, there are two $\E$-triangles
$$
R_C\rightarrow Q_C\rightarrow C\overset{}{\dashrightarrow}~~\textrm{and}~~
C\rightarrow R^C\rightarrow Q^C\overset{}{\dashrightarrow}
$$
satisfying $Q_C,Q^C\in \Q$ and $R_C,R^C\in \R$.
\end{itemize}
A cotorsion pair $(\Q,\R)$ is called \emph{hereditary} if  $\Q$ is closed under cocones of deflations and $\R$ is closed under cones of inflations.
\end{definition}

\begin{remark}
If $(\C,\E,\s)$ is an exact category, then the above cotorsion pair is just complete cotorsion pair in the sense of \cite{H2}.
\end{remark}

Now we state and prove our main result.

\begin{theorem}\label{main}
Assume that $(\mathcal{C}, \mathbb{E}, \mathfrak{s})$ is an extriangulated category satisfying {\rm Condition (WIC)}.
Let $(\Q, \widetilde{\R})$ and $(\widetilde{\Q}, \R)$ be two hereditary cotorsion pairs satisfying conditions $\widetilde{\R} \subseteq \R$, $\widetilde{\Q}\subseteq \Q$ and $\widetilde{\Q}\cap \R = \Q \cap \widetilde{R}$. Then there exists a unique thick class $\W$ for which $(\Q,\W,\R)$ is a Hovey triple. Moreover, this thick class $\W$ can be described in the two following ways:
\begin{align*}
   \W  &= \{\, X \in \C \, | \,  \text{there exists an $\E$-triangle} \, X \to R \to Q\dashrightarrow \, \text{ with} \, R \in \widetilde{\R } \, , Q \in \widetilde{\Q} \,\} \\
           &= \{\, X \in \C \, | \,  \text{there exists an $\E$-triangle } \, R' \to Q' \to X \dashrightarrow\, \text{ with} \, R' \in \widetilde{\R} \, , Q' \in \widetilde{\Q} \,\}.
          \end{align*}
\end{theorem}

\proof (1) We write
$$\W_1=\{\, X \in \C \, | \,  \text{there exists an $\E$-triangle } \, X \to R \to Q\dashrightarrow \, \text{ with} \, R \in \widetilde{\R } \, , Q \in \widetilde{\Q} \,\},$$
$$\W_2=\{\, X \in \C \, | \,  \text{there exists an $\E$-triangle } \, R' \to Q' \to X\dashrightarrow \, \text{ with} \, R' \in \widetilde{\R} \, , Q' \in \widetilde{\Q} \,\}.$$
Now show that $\W_1=\W_2$.

For any $X\in\W_1$, there exists an $\E$-triangle $X\to R\to Q\dashrightarrow$ where $R \in \widetilde{\R}$ and $Q \in \widetilde{\Q}$.
Since $(\Q,\widetilde{\R})$ is cotorsion pair, there exists an $\E$-triangle
$R'\to\Q'\to R\dashrightarrow$ where $R'\in\widetilde{\R}$ and $Q'\in\Q$.
By (ET4$\op$),  we have the following commutative diagram of $\E$-triangles.
$$\xymatrix{R'\ar[r]\ar@{=}[d] &M\ar[r]\ar[d] &X\ar@{-->}[r]\ar[d]&\\
R'\ar[r] & Q'\ar[r]\ar[d] & R\ar@{-->}[r]\ar[d] &\\
& Q\ar@{-->}[d]\ar@{=}[r] & Q\ar@{-->}[d]\\&&}$$
Since $\widetilde{\R}$ is closed under extensions, we have $Q'\in\widetilde{\R}$. It follows that
$Q'\in\Q \cap \widetilde{R}=\widetilde{\Q}\cap \R$ and then $Q'\in\widetilde{\Q}$.
Now because $Q',Q\in\widetilde{\Q}$ and $\widetilde{Q}$ is closed under cocones of deflations, we conclude that $M\in\widetilde{\Q}$.
Thus we find an $\E$-triangle
$R'\to M\to X\dashrightarrow$ where $R'\in\widetilde{\R}$ and $M\in\widetilde{\Q}$.
Hence $\W_1\subseteq\W_2$. A similar argument will show that $\W_2\subseteq\W_1$.

(2) We show that $\W$ is thick.

$\bullet$ We first prove that $\W$ is closed under direct summands.

Suppose $W\in\W$ and $X\xrightarrow{~i~}W\xrightarrow{~p~}X$ satisfies $pi=1_X$.
We want to show that $X\in\W$.

Since $W\in\W$, there exists an $\E$-triangle $W\overset{x}{\longrightarrow} \widetilde{R}\overset{y}{\longrightarrow} \widetilde{Q}\dashrightarrow$ where $\widetilde{R}\in\widetilde{\R}$ and $\widetilde{Q}\in\widetilde{\Q}$.
Since $(\widetilde{\Q}, \R)$ is a cotorsion pair, there exists an $\E$-triangle
$X\overset{k}{\longrightarrow} R\overset{h}{\longrightarrow} \widetilde{Q'}\dashrightarrow$ where $R\in\R$ and $\widetilde{Q'}\in\widetilde{\Q}$.

Applying the functor $\Hom_{\C}(-,\widetilde{R})$ to the $\E$-triangle $X\overset{k}{\longrightarrow} R\overset{h}{\longrightarrow} \widetilde{Q'}\dashrightarrow$, by Lemma \ref{lem}, we have the following exact sequence:
$$\Hom_{\C}(R,\widetilde{R})\xrightarrow{\Hom_{\C}(k,~\widetilde{R})}\Hom_{\C}(X,\widetilde{R})\xrightarrow{}\E(\widetilde{Q'},\widetilde{R})=0.$$
Thus there exists a morphism $j\colon R\to \widetilde{R}$ such that
$xi=jk$.

Applying the functor $\Hom_{\C}(-,R)$ to the $\E$-triangle $W\overset{x}{\longrightarrow} \widetilde{R}\overset{y}{\longrightarrow} \widetilde{Q}\dashrightarrow$, by Lemma \ref{lem}, we have the following exact sequence:
$$\Hom_{\C}(\widetilde{R},R)\xrightarrow{\Hom_{\C}(x,~R)}\Hom_{\C}(W,R)\xrightarrow{}\E(\widetilde{Q},R)=0.$$
Thus there exists a morphism $q\colon \widetilde{R}\to \R$ such that
$kp=qx$. By (ET3), we have the following commutative diagram
$$\xymatrix{
X \ar[r]^k \ar[d]^i& R\ar[r]^h \ar[d]^{j} & \widetilde{Q'}\ar@{-->}[r]\ar@{-->}[d]&\\
W \ar[r]^x \ar[d]^p & \widetilde{R}\ar[r]^y \ar[d]^{q} & \widetilde{Q}\ar@{-->}[r]\ar@{-->}[d]&\\
X\ar[r]^{k} & R \ar[r]^{h} & \widetilde{Q'}\ar@{-->}[r] &}$$
of $\E$-triangles.
It follows that $(1_R-qj)k=k-qjk=k-kpi=0$. So there exists a morphism
$t\colon \widetilde{Q'}\to R$ such that $th=1_R-qj$. That is to say, $1_R-qj$ factors through $\widetilde{Q'}$.
We claim that $1_R-qj$ factors through some object in $\widetilde{\Q}\cap \R = \Q \cap \widetilde{R}$.
In fact, Since $(\widetilde{\Q}, \R)$ is a cotorsion pair, there exists an $\E$-triangle
$\widetilde{Q'}\overset{u}{\longrightarrow} R'\overset{v}{\longrightarrow} \widetilde{Q''}\dashrightarrow$ where $R'\in\R$ and $\widetilde{Q''}\in\widetilde{\Q}$.
Since $\widetilde{\Q}$ is closed under extensions, we have $R'\in\widetilde{\Q}$ implies $R'\in\widetilde{\Q}\cap \R = \Q \cap \widetilde{R}$.
Hence $R'\in\widetilde{\R}$.

Applying the functor $\Hom_{\C}(-,R)$ to the $\E$-triangle $\widetilde{Q'}\overset{u}{\longrightarrow} R'\overset{v}{\longrightarrow} \widetilde{Q''}\dashrightarrow$, by Lemma \ref{lem}, we have the following exact sequence:
$$\Hom_{\C}(R',R)\xrightarrow{\Hom_{\C}(u,~R)}\Hom_{\C}(\widetilde{Q'},R)\xrightarrow{}\E(\widetilde{Q''},R)=0.$$
Thus there exists a morphism $\beta\colon R'\to R$ such that $t=\beta u$
and then $1_R-qj=th=\beta(uh)$.
It follows that $(q,\beta)\binom{j}{uh}=1_R$. Namely, the composition below is the identity $1_R$.
$$R\xrightarrow{~\binom{j}{uh}~}\widetilde{R}\oplus R'\xrightarrow{~(q,~\beta)~}R$$
This means that $R$ is a direct summand of $\widetilde{R}\oplus R'$. Since $\widetilde{\R}$ is closed under direct sums and direct
summands, we obtain $R\in\widetilde{\R}$.
This shows that $X\in\W$, as required.

$\bullet$ We show that $\W$ is closed under extensions.

Note that $\widetilde{\Q}\subseteq\W$ and $\widetilde{\R}\subseteq\W$, we have the following claim.

\textbf{Claim I}:  Let $R\overset{m}{\longrightarrow}Y\overset{n}{\longrightarrow}W\overset{}{\dashrightarrow}$ be an $\E$-triangle with $R\in\widetilde{\R}$ and $W\in\W$.
Then there exists a commutative diagram
 $$ \xymatrix{
  \widetilde{R} \ar@{}[dr]|{\circlearrowleft}\ar[d]\ar[r]  & \widetilde{R''}\ar@{}[dr]|{\circlearrowleft} \ar[d]\ar[r]  & \widetilde{R' } \ar[d]\ar@{-->}[r]  & \\
  \widetilde{Q }\ar@{}[dr]|{\circlearrowleft}\ar[d] \ar[r] & \widetilde{Q''}\ar@{}[dr]|{\circlearrowleft} \ar[d] \ar[r] & \widetilde{Q'} \ar[d] \ar@{-->}[r] &  \\
  R\ar@{-->}[d] \ar[r] & Y \ar@{-->}[d] \ar[r] & W\ar@{-->}[d] \ar@{-->}[r] &  \\
  &  & &  }
  $$
of $\E$-triangles, where $\widetilde{Q},\widetilde{Q'},Q''\in\widetilde{\Q}$ and $\widetilde{R},\widetilde{R'},R''\in\widetilde{\R}$.

Indeed, since $R,W\in\W$, there are two $\E$-triangles:
$$\widetilde{R}\overset{b}{\longrightarrow}\widetilde{Q}\overset{a}{\longrightarrow}R\overset{}{\dashrightarrow}
~\textrm{and}~\widetilde{R'}\overset{d}{\longrightarrow}\widetilde{Q'}\overset{c}{\longrightarrow}W\overset{}{\dashrightarrow}$$
where $\widetilde{R},\widetilde{R'}\in\widetilde{\R}$ and $\widetilde{Q},\widetilde{Q'}\in\widetilde{\Q}$.
Since $R\in\widetilde{\R}\subseteq\R$, we have $\E(\widetilde{Q'},R)=0$.

Applying the functor $\Hom_{\C}(\widetilde{Q'},-)$ to the $\E$-triangle $R\overset{m}{\longrightarrow}Y\overset{n}{\longrightarrow}W\overset{}{\dashrightarrow}$, by Lemma \ref{lem}, we have the following exact sequence:
$$\Hom_{\C}(\widetilde{Q'},Y)\xrightarrow{\Hom_{\C}(\widetilde{Q'},~n)}\Hom_{\C}(\widetilde{Q'},W)\xrightarrow{}\E(\widetilde{Q'},R)=0.$$
Thus there exists a morphism $w\colon \widetilde{Q'}\to Y$ such that $c=nw$.
Hence we have the following commutative diagram
$$\xymatrix@C=1.2cm{
\widetilde{Q} \ar[r]^{\binom{1}{0}\quad} \ar[d]^a & \widetilde{Q}\oplus \widetilde{Q'}\ar[r]^{\quad(0,~1)}\ar[d]^{(ma,~w)} & \widetilde{Q'}\ar@{-->}[r]\ar[d]^c&\\
R\ar[r]^{m} & Y \ar[r]^{n} & W\ar@{-->}[r] &}$$
of $\E$-triangles. Since $a,c$ are deflations, by Lemma \ref{y5}, there exists a deflation $b'\colon \widetilde{Q}\oplus \widetilde{Q'}\to Y$ which gives the following
morphism of $\E$-triangles.
$$\xymatrix@C=1.2cm{
\widetilde{Q} \ar[r]^{\binom{1}{0}\quad} \ar[d]^a & \widetilde{Q}\oplus \widetilde{Q'}\ar[r]^{\quad(0,~1)}\ar[d]^{b'} & \widetilde{Q'}\ar@{-->}[r]\ar[d]^c&\\
R\ar[r]^{m} & Y \ar[r]^{n} & W\ar@{-->}[r] &}$$
Since $b'$ is a deflation, there exists an $\E$-triangle $\widetilde{R''}\to\widetilde{Q}\oplus \widetilde{Q'}\to Y\dashrightarrow $.
By Lemma \ref{y8}, we have a commutative diagram
 $$ \xymatrix@C=1.2cm{
  \widetilde{R} \ar@{}[dr]|{\circlearrowleft}\ar[d]\ar[r]  & \widetilde{R''}\ar@{}[dr]|{\circlearrowleft} \ar[d]\ar[r]  & \widetilde{R' } \ar[d]\ar@{-->}[r]  & \\
  \widetilde{Q }\ar@{}[dr]|{\circlearrowleft}\ar[d] \ar[r]^{\binom{1}{0}\quad} & \widetilde{Q}\oplus \widetilde{Q'}\ar@{}[dr]|{\circlearrowleft} \ar[d] \ar[r]^{\quad(0,~1)} & \widetilde{Q'} \ar[d] \ar@{-->}[r] &  \\
  R\ar@{-->}[d] \ar[r] & Y \ar@{-->}[d] \ar[r] & W\ar@{-->}[d] \ar@{-->}[r] &  \\
  &  & &  }
  $$
of $\E$-triangles.
Since $\widetilde{\Q}$ is closed under direct sums and extensions, we have $\widetilde{Q}\oplus \widetilde{Q'}\in\widetilde{\Q}$.
Since $\widetilde{\R}$ is closed under extensions, we have $\widetilde{R''}\in\widetilde{\R}$.
So we have proved our claim.

Now we show that $\W$ is closed under extensions.
Let $$W\overset{}{\longrightarrow}Y\overset{}{\longrightarrow}W'\overset{}{\dashrightarrow}$$ be an $\E$-triangle with $W,W'\in\W$.
We need to show that $Y\in\W$ too.

Since $W\in\W$, there exists an $\E$-triangle $W\to R\to Q\dashrightarrow$ where $R\in\widetilde{\R}$ and $Q\in\widetilde{\Q}$.
By \cite[Proposition 3.15]{NP}, we have the following commutative diagram
$$\xymatrix{W\ar[r]\ar[d] &Y\ar[r]\ar[d] &W'\ar@{-->}[r]\ar@{=}[d]&\\
R\ar[r]\ar[d] & M\ar[r]\ar[d] & W'\ar@{-->}[r] &\\
 Q\ar@{-->}[d]\ar@{=}[r] & Q\ar@{-->}[d]&\\&&}$$
of $\E$-triangles.
By Claim I, there exists a commutative diagram
 $$ \xymatrix{
  \widetilde{R} \ar@{}[dr]|{\circlearrowleft}\ar[d]\ar[r]  & \widetilde{R''}\ar@{}[dr]|{\circlearrowleft} \ar[d]\ar[r]  & \widetilde{R' } \ar[d]\ar@{-->}[r]  & \\
  \widetilde{Q }\ar@{}[dr]|{\circlearrowleft}\ar[d] \ar[r] & \widetilde{Q''}\ar@{}[dr]|{\circlearrowleft} \ar[d] \ar[r] & \widetilde{Q'} \ar[d] \ar@{-->}[r] &  \\
  R\ar@{-->}[d] \ar[r] & M \ar@{-->}[d] \ar[r] & W\ar@{-->}[d] \ar@{-->}[r] &  \\
  &  & &  }
  $$
of $\E$-triangles, where $\widetilde{Q},\widetilde{Q'},Q''\in\widetilde{\Q}$ and $\widetilde{R},\widetilde{R'},R''\in\widetilde{\R}$.
By (ET4$\op$),  we have the following commutative diagram
$$\xymatrix{\widetilde{R''}\ar[r]\ar@{=}[d] &L\ar[r]\ar[d] &Y\ar@{-->}[r]\ar[d]&\\
\widetilde{R''}\ar[r] & \widetilde{Q''}\ar[r]\ar[d] & M\ar@{-->}[r]\ar[d] &\\
& Q\ar@{-->}[d]\ar@{=}[r] & Q\ar@{-->}[d]\\&&}$$
of $\E$-triangles.
Since $\widetilde{\Q}$ is closed under cocones of deflations, we have $L\in\widetilde{\Q}$.
Therefore $Y\in\W$. This show that $\W$ is closed under extensions.

$\bullet$ Now we show that $\W$ satisfies
the $2$-out-of-$3$ property.

Let $W\overset{}{\longrightarrow}W'\overset{}{\longrightarrow}Z\overset{}{\dashrightarrow}$ be an $\E$-triangle with $W,W'\in\W$.
We need to show that $Z\in\W$ too.  A dual argument shows that if $X\to W'\to W\dashrightarrow$ is an $\E$-triangle with $W,W'\in\W$, then $X\in\W$.

Since $W\in\W$, there exists an $\E$-triangle $W\to \widetilde{R}\to \widetilde{Q}\dashrightarrow$ where $\widetilde{R}\in\widetilde{\R}$ and $\widetilde{Q}\in\widetilde{\Q}$.
By \cite[Proposition 3.15]{NP}, we have the following commutative diagram
$$\xymatrix{W\ar[r]\ar[d] &W'\ar[r]\ar[d] &Z\ar@{-->}[r]\ar@{=}[d]&\\
\widetilde{R}\ar[r]\ar[d] & M\ar[r]\ar[d] & Z\ar@{-->}[r] &\\
 \widetilde{Q}\ar@{-->}[d]\ar@{=}[r] & \widetilde{Q}\ar@{-->}[d]&\\&&}$$
of $\E$-triangles.
Since $\W$ is closed under extensions, we have $M\in\W$.
So there exists an $\E$-triangle
$M\to \widetilde{R'}\to \widetilde{Q'}\dashrightarrow$ where $\widetilde{R'}\in\widetilde{\R}$ and $\widetilde{Q'}\in\widetilde{\Q}$.
By (ET4),
we have the following commutative diagram
$$\xymatrix{\widetilde{R}\ar[r]\ar@{=}[d] &M\ar[r]\ar[d] &Z\ar@{-->}[r]\ar[d]&\\
\widetilde{R}\ar[r] & \widetilde{R'}\ar[r]\ar[d] & L\ar@{-->}[r]\ar[d] &\\
& \widetilde{Q'}\ar@{-->}[d]\ar@{=}[r] & \widetilde{Q'}\ar@{-->}[d]\\&&}$$
of $\E$-triangles.
Since $\widetilde{\R}$ is closed under cones of inflations, we have $L\in\widetilde{\R}$.
Hence $Z\in\W$.
\medskip

This completes the proof that $\W$ is thick.

\medskip
(3) Now we show that $(\Q,\W,\R)$ is a Hovey triple. It suffices to show $\W\cap\R=\widetilde{R}$ and $\Q\cap\W=\widetilde{Q}$.

We only prove that $\W\cap\R=\widetilde{R}$.
Through similar assertions, we can prove $\Q\cap\W=\widetilde{Q}$.

It is clear that $\widetilde{\R}\subseteq\W\cap\R$.
Conversely, for any $X\in\W\cap\R$, there exists an $\E$-triangle
$X\to \widetilde{R}\to \widetilde{Q}\dashrightarrow$ where $\widetilde{R}\in\widetilde{\R}$ and $\widetilde{Q}\in\widetilde{\Q}$.
Since $X\in\R$, we have $\E(\widetilde{Q},X)=0$. We conclude that the above $\E$-triangle splits.
Thus $X$ is a direct summand of $\widetilde{R}$. It follows that $X\in\widetilde{\R}$.
This completes the proof that $(\Q,\W,\R)$ is a Hovey triple.
\medskip

(4) Now we show that $\W$ is a unique.

Suppose that $(\Q, \V, \R)$ is any other Hovey triple, we must have $\V\cap\R =\widetilde{\R}= \W\cap\R$
and $\Q\cap\V = \widetilde{Q}=\Q\cap\W$. By  the definitions of $\W$ and $\V$, we have $\W= \V$.   \qed

\medskip

When our main result is applied to an exact category, we have the following.

\begin{corollary}\emph{\cite[Theorem 1.2]{G2}}
Assume that $\C$ is a weakly idempotent complete exact category.
Let $(\Q, \widetilde{\R})$ and $(\widetilde{\Q}, \R)$ be two hereditary cotorsion pairs satisfying conditions $\widetilde{\R} \subseteq \R$, $\widetilde{\Q}\subseteq \Q$ and $\widetilde{\Q}\cap \R = \Q \cap \widetilde{R}$. Then there exists a unique thick class $\W$ for which $(\Q,\W,\R)$ is a Hovey triple. Moreover, this thick class $\W$ can be described in the two following ways:
\begin{align*}
   \W  &= \{\, X \in \C \, | \,  \text{there exists an exact sequence} \, X \rightarrowtail R \twoheadrightarrow Q \, \text{ with} \, R \in \widetilde{\R } \, , Q \in \widetilde{\Q} \,\} \\
           &= \{\, X \in \C \, | \,  \text{there exists an exact sequence } \, R' \rightarrowtail Q' \twoheadrightarrow X\, \text{ with} \, R' \in \widetilde{\R} \, , Q' \in \widetilde{\Q} \,\}.
          \end{align*}
\end{corollary}

When our main result is applied to a triangulated category, we have the following.

\begin{corollary}
Assume that $\C$ is a triangulated category.
Let $(\Q, \widetilde{\R})$ and $(\widetilde{\Q}, \R)$ be two hereditary cotorsion pairs satisfying conditions $\widetilde{\R} \subseteq \R$, $\widetilde{\Q}\subseteq \Q$ and $\widetilde{\Q}\cap \R = \Q \cap \widetilde{R}$. Then there exists a unique thick class $\W$ for which $(\Q,\W,\R)$ is a Hovey triple. Moreover, this thick class $\W$ can be described in the two following ways:
\begin{align*}
   \W  &= \{\, X \in \C \, | \,  \text{there exists a triangle} \, X\to R \to Q \to X[1]\, \text{ with} \, R \in \widetilde{\R } \, , Q \in \widetilde{\Q} \,\} \\
           &= \{\, X \in \C \, | \,  \text{there exists a triangle} \, R' \to Q' \to X\to R'[1]\, \text{ with} \, R' \in \widetilde{\R} \, , Q' \in \widetilde{\Q} \,\}.
          \end{align*}
\end{corollary}

\textbf{Panyue Zhou}\\
College of Mathematics, Hunan Institute of Science and Technology, 414006, Yueyang, Hunan, P. R. China.\\
and\\
D\'{e}partement de Math\'{e}matiques, Universit\'{e} de Sherbrooke, Sherbrooke,
Qu\'{e}bec J1K 2R1, Canada.\\
E-mail: \textsf{panyuezhou@163.com}

\end{document}